\documentclass[11pt,a4paper]{article}

\usepackage[margin=1in]{geometry}
\usepackage[T1]{fontenc}
\usepackage[utf8]{inputenc}
\usepackage{lmodern}
\usepackage{microtype}

\usepackage{amsmath,amssymb}
\usepackage{adjustbox}
\usepackage{booktabs}
\usepackage[short,nocomma]{optidef}
\usepackage{graphicx}
\usepackage[dvipsnames]{xcolor}
\usepackage[colorlinks=true,linkcolor=MidnightBlue,citecolor=MidnightBlue,urlcolor=MidnightBlue]{hyperref}
\usepackage[capitalize,nameinlink]{cleveref}

\makeatletter
\def\input@path{{./}{sections/}}
\makeatother
\graphicspath{{./}}

\newcommand{\nextvar}{\texttt{NextVar}}
\newcommand{\vehiclevar}{\texttt{VehicleVar}}
\newcommand{\activevar}{\texttt{ActiveVar}}
\newcommand{\timevar}{\texttt{Time}}
\newcommand{\loadvar}{\texttt{Load}}
\newcommand{\true}{\texttt{True}}
\newcommand{\false}{\texttt{False}}

\title{Paratransit Optimization with Constraint Programming:\\A Case Study in Savannah, Georgia}

\author{
  Liam Jagrowski \quad Kevin Dalmeijer \quad Tinghan Ye \quad Pascal Van Hentenryck \\[4pt]
  \normalsize H. Milton Stewart School of Industrial and Systems Engineering, Georgia Tech, USA \\
  \normalsize \texttt{\{ljagrowski3, dalmeijer, joe.ye, pvh\}@gatech.edu}
}

\date{}

\begin{document}

\maketitle

\begin{abstract}
Paratransit services are vital for individuals who cannot use fixed-route public transit, including those with disabilities.
Optimizing these services is essential for transit agencies to deliver high-quality service efficiently.
This paper introduces a Constraint Programming (CP) model to jointly optimize route planning and shift scheduling for paratransit operations, along with practical guidance for real-world implementation.
A case study in Savannah, Georgia, demonstrates that the new approach is competitive with a recently proposed, highly effective AI-accelerated column generation framework, and significantly increases the number of requests served compared to current practices.
The method is also easier to implement and provides an inherently practical solution for transportation planners.
CP further provides the flexibility to optimize schedules without requiring shifts to start exactly on the hour, yielding an additional 5\% improvement in the number of requests served.
\end{abstract}

\section{Introduction}
\label{sec:introduction}

Paratransit, also known as Intermediate Public Transport or Community Transport, is designed to supplement the public transit system with flexible and individualized rides.
These services are critical to provide mobility for the elderly, disabled, and other individuals who would have challenges using conventional public transit.
Paratransit systems may be performed with a variety of vehicles, ranging from small sedans to large converted vans \cite{fuparatransit}.
Riders call in or use an app to request door-to-door service at a given time, after which paratransit planners create routes to serve as many requests as possible.

Paratransit services bring a clear value to the community and laws may be in place to make sure that they are provided.
In the US, for example, the \textit{Americans with Disabilities Act} compels transit authorities that provide fixed-route options to also offer a paratransit service with a level of service comparable to that of the fixed-route system \cite{dotada}.
However, providing this level of service is often challenging due to the growing demand from an aging population and limited operational resources.

In the literature, there has been significant efforts to optimize paratransit systems to navigate these challenges \cite{fuparatransit,lu2024boosting,zhang21pt}.
For example, \cite{fuparatransit} optimize the size and composition of the paratransit fleet, and \cite{zhang21pt} provide mixed integer programming models to minimize user waiting time as well as operating costs.
In terms of modeling, paratransit services differ significantly from other forms of transportation.
Each route is customized, and there is often no fixed timetable.
This unique nature of paratransit requires distinct approaches and solutions.

\cite{lu2024boosting} is a seminal work that first introduces the problem of jointly optimizing rider trip planning and crew scheduling with realistic constraints for the paratransit setting considered here.
Rather than assuming that the driver shifts are given, they are optimized to best align with the demand.
Furthermore, the authors support travel plans that consist of multiple trips, e.g., to and from the hospital.
This problem is solved with an algorithm that combines column generation with machine learning to obtain near optimal solutions.
While this integration of artificial intelligence and operations research is effective and significantly outperforms current practice, it is also rather complex and relies on historical data to train the underlying machine learning model.
Transit agencies, however, also need optimization tools that are maintainable, flexible, and resilient to changes in operating policy and time resolution.

To complement the work by \cite{lu2024boosting}, this paper introduces a Constraint Programming (CP) model that is easier to implement and that provides near-identical performance without relying on historical data.
Practical guidance is provided to support the implementation of the new model in real-world systems, positioning it as a practical bridge between advanced optimization and on-the-ground deployment.
As an additional benefit, CP provides the flexibility to optimize schedules without requiring shifts to start exactly on the hour, making the approach more robust to changes in scheduling granularity.
Enabling this option allows for a 5\% increase in requests served, and the new CP model outperforms \cite{lu2024boosting} in this setting.

\subsection{Contributions}
The contributions of this paper can be summarized as follows:
\begin{enumerate}
    \item The paper introduces a CP model to jointly optimize route planning and shift scheduling for paratransit operations, along with practical guidance for real-world implementation.
    \item A case study in Savannah demonstrates that the new approach is competitive with a recently proposed, highly effective AI-accelerated column generation framework, and significantly increases the number of requests served compared to current practices.
	\item The new model makes it possible to relax the requirement that shifts start exactly on the hour, resulting in an additional 5\% increase in the number of requests served, an improvement that was previously unattainable.
\end{enumerate}

\noindent The remainder of this paper is organized as follows.
 Section~\ref{sec:literature} provides a literature review.
 Section~\ref{sec:methodology} introduces the model and discusses the implementation details.
 The case study in Savannah is introduced in Section~\ref{sec:casestudy}, and Section~\ref{sec:results} provides the results.
 Section~\ref{sec:conclusion} ends the paper with a conclusion and directions for future research.

\section{Literature Review}
\label{sec:literature}

This section reviews the relevant literature that frames the research presented. The review first discusses the paratransit optimization problem as a complex variant of the Dial-a-Ride Problem (DARP). It then examines the literature on the joint optimization of vehicle routes and crew schedules. Finally, the study is situated within the context of practical implementation using modern open-source solvers like Google OR-Tools.

\subsection{Paratransit Routing as a Dial-a-Ride Problem}
The Dial-a-Ride Problem (DARP) and the more general Pickup-and-Delivery Problem (PDP) provide the mathematical foundation for paratransit services. The goal is to design minimum-cost vehicle routes to serve a set of user requests, subject to time windows and vehicle capacity. This problem class is well-studied, with exact methods based on mixed-integer programming (MIP) providing optimal solutions for small to medium-sized instances \cite{cordeau2006branch,rist2021new}. In parallel, CP has emerged as a powerful paradigm for handling the rich, non-linear, and logical side constraints that are often difficult to express in classical MIP formulations \cite{kilby2006vehicle,lam2016branch,lam2020joint}. Specifically for the dial-a-ride setting, \cite{cappart2018constraint} develop a CP model for patient transportation problems, and \cite{thomas2020insertion} introduce insertion sequence variables for hybrid routing and scheduling. The present work builds on the same modeling tradition but targets an integrated paratransit setting using the OR-Tools routing library, focusing on a practical implementation that supports both fixed and flexible shifts.

However, real-world paratransit operations introduce domain-specific complexities that extend beyond classical DARP models. These include diverse rider needs, regulatory requirements, and a strong emphasis on service quality metrics like punctuality and minimal ride times \cite{paquette2009quality}. This has led to the development of sophisticated decomposition techniques and specialized metaheuristics to tackle large-scale, practical instances \cite{karabuk2009nested,kirchler2013granular}. As surveyed by \cite{ho2018survey} and \cite{molenbruch2017typology}, the field has produced a wide array of exact, heuristic, and hybrid algorithms tailored to these challenges. A critical recent trend is bridging the gap between theory and practice, exemplified by work like \cite{pavia2024deploying}, who successfully deployed an optimization-based scheduler with a US transit agency, demonstrating field-tested improvements.

\subsection{Joint Optimization of Vehicle Routing and Crew Scheduling}
While vehicle routing and driver shift planning are often addressed sequentially in practice, a growing body of literature demonstrates that integrated optimization yields superior results. Seminal work by \cite{huisman2005multiple} showed that combining vehicle and crew scheduling can significantly reduce overall costs and resource requirements. More recently, research has focused on enhancing these integrated models to handle real-world uncertainty and disruptions \cite{amberg2023robust}.

Despite these advantages, integrated models remain rare in the paratransit context due to the immense computational complexity. A notable recent exception is \cite{lu2024boosting}, who present an integrated model for joint trip and shift scheduling using a sophisticated column generation method enhanced by a graph neural network. While powerful, the complexity of such methods can pose a barrier to implementation for many transit agencies, highlighting a clear need for models that are both comprehensive and practical.

\subsection{Practical Implementation with Google OR-Tools}
The present study aims to develop a model that is not only effective but also practical to implement. Google's OR-Tools suite has become a powerful and accessible tool for this purpose. While its use in academic paratransit literature is still emerging, its effectiveness is well-documented in analogous routing problems \cite{silva2023capacitated}.
For example, \cite{alves2021solving} used it for real-world waste collection -- a capacitated pickup problem -- highlighting the power of its built-in metaheuristics. Other studies confirm that OR-Tools can rapidly generate high-quality solutions for various VRPs \cite{khanna2025pdptw,silva2025hybrid,silva2023capacitated}.

The work most relevant to this study is by \cite{pavia2024deploying}, which details a successful application of OR-Tools in a real-world paratransit setting. This research builds directly on that work by using the same underlying solver technology. However, a key gap is addressed by extending the problem formulation to jointly optimize vehicle routes and crew schedules. The model developed in this paper accomplishes this using CP to manage the combined complexity while introducing novel flexibility, such as allowing unconstrained shift start times, to further improve operational efficiency.

\section{Methodology}
\label{sec:methodology}

This section formally describes the problem of optimizing paratransit operations and models it as a CP problem.
It also documents the implementation choices used to solve this CP through the Google OR-Tools Routing Library in the experiments.
These choices are presented as part of the experimental setup rather than as universal best practices, and Section~\ref{sec:results} empirically evaluates the impact of the proposed acceleration technique.

The problem solved here matches the Joint Rider Trip Planning and Crew Shift Scheduling Problem (JRTPCSSP) as introduced by \cite{lu2024boosting}.
Following prior work and consistent with the current workflow of the partner agency, the proposed model assumes a reservation-based (static) setting in which requests are received one day in advance, and the day-ahead schedule is optimized with complete information.
However, the CP model presented here can be readily adapted to a dynamic setting by employing a rolling-horizon framework, which represents an interesting direction for future research.

\subsection{Problem Description}
The goal of the JRTPCSSP is to serve a set of passenger requests $R=\{1, \hdots, n\}$ with a set of vehicles $K = \{1, \hdots, m\}$, simultaneously optimizing vehicle routes and the shifts during which the vehicles are active.
Let $G=(V,A)$ be a graph with vertices $V$ and directed arcs $A$.
Each request $i \in R$ is associated with a pickup node $i \in V$ and a drop-off node $n + i \in V$.
Furthermore, $V$ contains starting depot node $2n+k \in V$ for each vehicle $k\in K$ and an ending depot node $2n+m+k \in V$ at which each vehicle $k$ ends.
For convenience, let $P=\{1,\hdots, n\}$ be the set of pickup nodes, let $D=\{n+1,\hdots, 2n\}$ be the set of drop-off nodes, let $S = \{2n+1,\hdots, 2n+m\}$ be the set of starting depots, and let $T = \{2n+m+1, \hdots, 2n+2m\}$ be the set of ending depots.
It follows that $V = S \cup P \cup D \cup T$.
Furthermore, let $V' = P \cup D$ be the set of all request nodes.
Directed arcs $(i, j) \in A$ indicate travel from node $i \in V$ to node $j \in V$.
Arcs are defined out of the starting depots ($S$ to $V'$), between the request nodes ($V'$ to $V'$), and into the ending depots ($V'$ to $T$).
Finally, arcs $(2n + k, 2n + m + k)$ $\forall k \in K$ are defined between the starting depots and the corresponding ending depots to indicate when vehicles are inactive and stay at the depot.

The route of vehicle $k \in K$ is represented as a simple path in $G$ from starting depot $2n+k$ to ending depot $2n+m+k$ that satisfies additional constraints, such as time constraints and capacity constraints.
With regards to time, each node $i \in V$ is associated with a service time $s_i \ge 0$ for boarding/unboarding, and each arc $(i,j) \in A$ is associated with a travel time $t_{ij} \ge 0$ to drive between locations. Two distinct nodes $i$ and $j$ may correspond to the same physical location, in which case $t_{ij}$ is zero and any boarding/unboarding interaction is captured by the service times.
Throughout this paper it is further assumed that all input data is integer.
This is not restrictive, as the problem can be rescaled by choosing an appropriate discretization of time, e.g., 1-minute versus 5-minute intervals.
Services times and travel times accumulate over the route, and each node $i\in V$ that is visited must start service within a prespecified time window $[a_i, b_i]$.
Vehicles are allowed to wait until the time window opens, but late arrivals are forbidden.
Furthermore, vehicles are restricted to leave the starting depot at one of the candidate starting times in the set $\mathcal{T}$, e.g., shifts only start at the top of the hour.
Finally, the maximum number of working hours per vehicle, i.e., the time between leaving the starting depot in $S$ and arriving at the ending depot in $T$, is limited by the maximum shift duration $L \ge 0$.

Routes are also required to satisfy pickup-and-delivery and capacity constraints.
Each pickup $i \in P$ must be followed by a drop-off $n + i \in D$ later in the route, and drop-off nodes cannot be visited without the corresponding pickup first.
Note that passengers do not need to be dropped off immediately.
For example, a vehicle can make multiple pickups followed by multiple drop-offs.
Each node $i \in V$ is associated with a demand $d_i$ that is positive when passengers enter the vehicle ($i \in P$) and negative when they exit ($i \in D$). In the experiments, each pickup demand is exactly canceled by its drop-off, i.e., $d_{n+i} = -d_i$ for all $i \in R$, so the load of every vehicle returns to zero at its terminal depot.
At all times, the number of people in each vehicle must remain below the vehicle capacity $Q \ge 0$.
Between different routes, each request can be served by at most one vehicle.

As it may be impossible to serve all requests on a given day, the objective is set to minimize the number of \emph{unserved requests}.
\cite{lu2024boosting} make the important observation that passengers often make multiple connected requests, for example a trip to the hospital and a return trip home.
To avoid that people get stranded, a constraint is added to serve either all or none of the connected requests.
That is, let $U$ be the set of requesters, and let $R_u \subseteq R$ be the set of requests made by $u \in U$ (such that the sets $R_u$ partition $R$).
Then either all requests in $R_u$ are served, or a penalty of $\lvert R_u \rvert$ is incurred: one for each unserved request.
Note that connected requests are allowed to be served by different vehicles, and these constraints thus span multiple routes.

In conclusion, a solution to the JRTPCSSP is represented by a set of paths in $G$ that satisfy the constraints above, as well as a set of times at which each node is visited. The paths define the routes of the vehicles and the trips of the passengers, while the departure and arrival times at the depots define the shifts for the drivers.

\subsection{Model}
The model in Figure~\ref{fig:cp_model} presents a CP formulation of the JRTPCSSP.
The model uses a set of decision variables that are commonly used in the CP literature to solve vehicle routing problems \cite{ortools_routing,kilby2006vehicle,lam2016branch}.
In particular, the notation in this paper is inspired by \cite{ortools_routing}.
For each node $i \in S \cup V'$, the successor variable $\nextvar(i) \in V' \cup T$ points to the next node that will be visited by the same vehicle, or it points to itself if the node is not visited at all.
Whether nodes are visited is also captured by the $\activevar(i) \in \{\true, ~\false\}$ variables for all $i \in S \cup V'$.
Finally, each node is associated with a vehicle $k \in K$ or with a non-existing vehicle 0 in case the node is inactive.
This is captured by the $\vehiclevar(i) \in \{0\} \cup K$ variables for all $i \in V$.

Constraints \eqref{eq:model:nextvar}-\eqref{eq:model:activevehicle} are standard vehicle routing constraints to ensure that the \nextvar{}, \activevar{}, and \vehiclevar{} variables represent valid routes.
A brief presentation is included here for completeness, and the reader is referred to \cite{kilby2006vehicle} for a detailed discussion.
It is also worth noting that Constraints \eqref{eq:model:nextvar}-\eqref{eq:model:activevehicle} may be omitted if they are already handled implicitly by the solver.
This is the case for the Google OR-Tools Routing Library \cite{ortools_routing}, for example, which is used for the experiments in this paper.
Equation~\eqref{eq:model:nextvar} defines the \nextvar{} successor variables and Constraint~\eqref{eq:model:alldifferent} forces all successors to be unique.
Equation~\eqref{eq:model:activevar} defines the \activevar{} variables, and Constraints~\eqref{eq:model:activenext} ensure that a node is marked active if and only if it has a successor other than itself.
Constraint \eqref{eq:model:nocycle} prevents cycles between active request nodes, while still allowing self-loops created by inactive nodes (see \cite{kilby2006vehicle} for more information about the \texttt{NoCycle} constraint).
The \vehiclevar{} variables are defined by Equation~\eqref{eq:model:vehiclevar}.
Constraints~\eqref{eq:model:depotvehicles} ensure that the depots are associated with the correct vehicles, and Constraints~\eqref{eq:model:nextvehicle} state that subsequent visits belong to the same vehicle.
Finally, Constraints~\eqref{eq:model:activevehicle} set \vehiclevar{} to zero for inactive nodes.

The remainder of the model in Figure~\ref{fig:cp_model} is specific to the JRTPCSSP.
The objective is to minimize the number of unserved requests.
If any of the connected requests $R_u$ of passenger $u \in U$ fail to be served, then all requests of this passenger are penalized for a total penalty of $\lvert R_u \rvert$. Formally, the indicator $I_{(\activevar(i)=\false \textrm{ for any } i \in R_u)}$ in Objective~\eqref{eq:model:obj} equals $1$ if there exists at least one $i \in R_u$ with $\activevar(i)=\false$ and equals $0$ when all requests in $R_u$ are active.
Time and capacity are modeled with two additional sets of variables: variables $\timevar(i)$ represent the time at which node $i \in V$ is served, and variables $\loadvar(i)$ represent the load of the vehicle on arrival at node $i \in V$.
Constraints~\eqref{eq:model:progress_time} are element constraints that define how time progresses for active nodes.
Equation~\eqref{eq:model:time_windows} enforces the time windows, and Equation~\eqref{eq:model:start_times} states that vehicles can only start at one of the candidate times in $\mathcal{T}$.
Furthermore, Constraints~\eqref{eq:model:duration} enforce the maximum shift length $L$ by limiting the time between leaving a starting depot and arriving at the corresponding ending depot.
The pickup and drop-off structure is imposed by Constraints~\eqref{eq:model:pd_vehicle} and \eqref{eq:model:pd_time}.
The former states that matching pickup and drop-off nodes are served by the same vehicle, or both remain unserved, while the latter forces the pickup to take place before the drop-off.
The next three constraints deal with vehicle capacity: element constraints~\eqref{eq:model:progress_load} define how the load of the vehicle is updated (recall that $d_i < 0$ for drop-offs), Equation~\eqref{eq:model:capacity} limits the load to the vehicle capacity $Q$, and Equation~\eqref{eq:model:load_start} defines the starting loads to be zero.

\begin{figure}[t]
\centering
\begin{mini!}
%
    % Variable
    {}
%
    % Objective
    {\sum_{u \in U} \lvert R_u \rvert I_{\left(\activevar(i)=\false \textrm{ for any } i \in R_u\right)}, \label{eq:model:obj}}
%
    % Model label
    {\label{model}}
%
    % Function closure for objective
    {}
%
    % Constraints
%
	\addConstraint
	{}
	{\nextvar{(i)} \in V' \cup T\label{eq:model:nextvar}}
	{\forall i \in S \cup V',}
	\addConstraint
	{}
	{\texttt{AllDifferent}\left(\nextvar{(i)}; i \in S \cup V'\right), \label{eq:model:alldifferent}}
	{}
	\addConstraint
	{}
	{\activevar{(i)} \in \{\texttt{True}, \texttt{False}\} \label{eq:model:activevar}}
	{\forall i \in S \cup V',}
	\addConstraint
	{}
	{\activevar{(i)} \Longleftrightarrow \left(\nextvar{(i)} \neq i \right) \label{eq:model:activenext}}
	{\forall i \in S \cup V',}
	\addConstraint
	{}
	{\texttt{NoCycle}\left(\nextvar{(i)}; i \in V' \vert \activevar{(i)}\right), \label{eq:model:nocycle}}
	{}
	\addConstraint
	{}
	{\vehiclevar{(i)} \in \{0\} \cup K \label{eq:model:vehiclevar}}
	{\forall i \in V,}
	\addConstraint
	{}
	{\vehiclevar{(2n+k)} = \vehiclevar{(2n + m + k)} = k \qquad \label{eq:model:depotvehicles}}
	{\forall k \in K,}
	\addConstraint
	{}
	{\vehiclevar{(\nextvar{(i)})} = \vehiclevar({i}) \label{eq:model:nextvehicle}}
	{\forall i \in S \cup V',}
	\addConstraint
	{}
	{\activevar{(i)} \Longleftrightarrow \left(\vehiclevar{(i)} \neq 0\right) \label{eq:model:activevehicle}}
	{\forall i \in S \cup V',}
    \addConstraint
    {}
    {\activevar{(i)} \Longrightarrow\notag}
    {}
    \addConstraint
	{}
	{\qquad \timevar(\nextvar(i)) \ge \timevar(i) + s_i + t_{i, \nextvar(i)} \label{eq:model:progress_time}}
	{\forall i \in S \cup V',}
    \addConstraint
    {}
    {\timevar(i) \in [a_i, b_i] \label{eq:model:time_windows}}
    {\forall i\in V,}
    \addConstraint
    {}
    {\timevar(i) \in \mathcal{T} \label{eq:model:start_times}}
    {\forall i \in S,}
    \addConstraint
    {}
    {\timevar(2n + m + k) - \timevar(2n + k) \le L \label{eq:model:duration}}
    {\forall k \in K,}
    \addConstraint
    {}
    {\vehiclevar(i) = \vehiclevar(n+i) \label{eq:model:pd_vehicle}}
    {\forall i\in R,}
    \addConstraint
    {}
    {\timevar(i) \le \timevar(n+i) \label{eq:model:pd_time}}
    {\forall i\in R,}
    \addConstraint
    {}
    {\activevar{(i)} \Longrightarrow \loadvar(\nextvar(i)) = \loadvar(i) + d_i \label{eq:model:progress_load}}
    {\forall i \in S \cup V',}
    \addConstraint
    {}
    {\loadvar(i) \in [0, Q] \label{eq:model:capacity}}
    {\forall i\in V,}
    \addConstraint
    {}
    {\loadvar(i) = 0 \label{eq:model:load_start}}
    {\forall i \in S.}
\end{mini!}%
\caption{CP Model for the JRTPCSSP.}
\label{fig:cp_model}
\end{figure}

\subsection{Implementation Details}
While the model presented above is general, practical performance still depends on how it is instantiated in a specific solver.
To support reproducibility, this paper documents the implementation used to solve the model in Figure~\ref{fig:cp_model} through the Google OR-Tools Routing Library \cite{ortools_routing}.
This library was chosen because it is open source, easy to use, and built on a CP solver.
The library handles standard vehicle routing constraints  \eqref{eq:model:nextvar}-\eqref{eq:model:activevehicle} by default and allows users to model various vehicle routing variants.
However, some of these implementation choices can have a substantial effect on performance in this setting.

The routing library exposes the underlying CP solver and allows users to directly add a variety of constraints.
However, it was found in preliminary experiments that adding constraints to the CP solver directly can significantly hurt performance, and in this setting it proved preferable to leverage the functions provided specifically by the routing library.
A good example is the implementation of the Objective~\eqref{eq:model:obj}.
\cite{lu2024boosting} support this objective by adding constraints that force the requests in $R_u$ to either all be served or all be unserved, and then penalize individual pickup nodes that remain unserved.
In the current notation, this corresponds to the constraints
% \\
\begin{equation}
    \label{eq:request_sets_naive}
    \activevar(i) = \activevar(j) \quad \forall u \in U, i, j \in R_u, i < j.
\end{equation}
% \\
In our preliminary tests, adding Constraints~\eqref{eq:request_sets_naive} to the CP solver directly led to poor performance, presumably because it interferes with local search.
With the above constraints it is only feasible to add a request to a route if all connected requests are added at the same time.
As an example, this causes the neighborhood operator \href{https://github.com/google/or-tools/blob/stable/ortools/constraint_solver/routing_neighborhoods.h}{\texttt{MakePairActiveOperator}} to fail for all $\lvert R_u \rvert > 2$, as it only inserts two currently-unserved requests at a time.
In the implementation used here, this issue is avoided by using the routing library's function \texttt{AddDisjunction()} instead.
For each connected set $R_u$, a single disjunction is posted over the pickup nodes in $R_u$ that allows up to $\lvert R_u \rvert$ of them to be active and that contributes a penalty of $\lvert R_u \rvert$ to the objective whenever fewer than $\lvert R_u \rvert$ are active. This recovers the indicator term in Objective~\eqref{eq:model:obj}: a partial group incurs the full primary penalty, exactly as if no member were served.\footnote{At the time of writing, the behavior of \texttt{AddDisjunction()} does not match the official documentation. The authors have verified that the actual behavior matches this paper and have contributed to the \href{https://github.com/google/or-tools/issues/1348}{issue on GitHub}.}
If a solution is returned with any request sets $R_u$ that are partially served, these requests can simply be removed in post-processing without changing the objective value.
Hence, this is without loss of optimality.

The time constraints \eqref{eq:model:progress_time}-\eqref{eq:model:time_windows} and capacity constraints \eqref{eq:model:progress_load}-\eqref{eq:model:load_start} are implemented as two \emph{dimensions}.
Dimensions in Google OR-Tools are similar to \emph{resources} in the vehicle routing literature (e.g., see \cite{IrnichDesaulniers2005-ShortestPathProblems}).
Each dimension defines how the resource is accumulated and defines the resource bounds at each node.
Constraints~\eqref{eq:model:start_times} are implemented by directly restricting the allowed values of $\timevar(i)$ through \texttt{SetValues()} and Constraints~\eqref{eq:model:duration} are supported through \texttt{Set\-Span\-Upper\-Bound\-For\-Vehicle()}.
Pickup and drop-off constraints~\eqref{eq:model:pd_vehicle}-\eqref{eq:model:pd_time} are added directly to solver.
In the implementation used here, the solver is also notified of this structure through \texttt{AddPickupAndDelivery()} to allow for better performance.
Finally, all nodes are mandatory to visit by default.
\texttt{AddDisjunction()} already makes the pickup nodes optional, and single-node disjunctions are also added for the drop-off nodes.
No additional penalty is used, as the penalties are already captured by the pickup nodes.
Finally, the local search metaheuristic \texttt{GENERIC\_TABU\_SEARCH} is enabled to apply tabu search \cite{Glover1989-TabuSearch—PartI} on the objective value of the solution to escape local minima. Section~\ref{sec:sensitivity} provides empirical support for these choices by comparing alternative metaheuristics and first-solution strategies.

\subsection{Acceleration Technique} \label{sec:acceleration}
From a computational perspective, one of the challenges of minimizing the number of unserved requests is that the solver has less incentive to optimize the \emph{order} of the requests on each route.
After all, successfully finding a faster route that serves the same requests does not improve the objective value.
However, optimizing the use of time in the current routes can greatly improve the search by making it easier to insert additional requests in the future.
Section~\ref{sec:results} evaluates this modeling choice empirically.

To take time into account during the solve, this paper proposes to include the total service and travel time as a secondary objective.
That is, Objective \eqref{eq:model:obj} is replaced by
% \\
\begin{equation}
\label{eq:obj:better}
    \hspace{2em} \sum_{\mathclap{\substack{i \in V' \cup S\\j = \nextvar(i)\\i\neq j}}} \quad\left(s_i + t_{ij}\right) + M \sum_{u \in U} \lvert R_u \rvert I_{\left(\activevar(i)=\false \textrm{ for any } i \in R_u\right)}.
\end{equation}
% \\
The constant $M$ is chosen sufficiently large such that serving requests is prioritized and the optimization problem is equivalent to the original \cite{sherali1982equivalent}.
In Google OR-Tools, the first term in Equation~\eqref{eq:obj:better} is easily implemented through the function \texttt{SetArcCostEvaluatorOfAll}\allowbreak\texttt{Vehicles()}.

\section{Case Study}
\label{sec:casestudy}

Chatham Area Transit (CAT) is the transit authority serving the Savannah metropolitan area. CAT offers fixed-route bus transit as well as paratransit and ferry service.
The paratransit service, known as \emph{CAT Mobility}, is crucial to its riders. This is evidenced by a relatively small decrease in usage during the pandemic: while fixed-route and ferry ridership decreased by 35\% and 42\%, respectively, from 2020 to 2021, paratransit services saw a much smaller 16\% drop \cite{CATplan}.
At the same time, CAT Mobility is becoming more expensive to operate, with a 36\% increase in operating cost per passenger trip between 2017 and 2021 \cite{CATState}.

These numbers imply a clear urgency to optimize CAT's paratransit operations to maximize the level of service with the limited resources available.
To evaluate the performance of the new CP model, this paper uses a real-world dataset derived from CAT Mobility that spans three pre-Christmas workweeks in December 2019.
These are the same instances used by \cite{lu2024boosting}, reused here deliberately to enable a fair comparison; the experimental contributions of this paper are the CP implementation, the acceleration ablation, and the flexible-shift analysis, rather than a new dataset.
Multi-trip passengers are quite frequent in this dataset, so connected requests arise regularly and introduce substantial inter-route complexity.
During these three weeks, CAT was only able to serve 81\% of requests.
This number can likely be attributed to a variety of factors.
This includes the limitations of the current planning tool to address this challenging optimization problem.
As reported by \cite{lu2024boosting}, CAT's existing operation relies on a commercial heuristic scheduler with predetermined driver shifts, and producing a daily plan typically takes hours; throughout this paper, the ``CAT'' baseline refers to the operational schedules actually produced by that system on the studied days.
Another contributing factor may be a shortage of drivers, as mentioned in the \emph{State of the System} report published by CAT in 2023 \cite{CATState}.
In any case, this provides a strong motivation to optimize the system.

\subsection{Experimental Settings}
The model in Figure~\ref{fig:cp_model} is implemented with Google OR-Tools v9.10.4067 in Python 3.9.
Experiments are conducted on a Linux machine with dual Intel Xeon Gold 6226 CPUs on the PACE Phoenix cluster \cite{PACE2017-PartnershipAdvancedComputing},
using up to eight cores and 64GB of RAM and a time limit of 30 minutes per instance.
These settings match \cite{lu2024boosting} to allow for a fair comparison.

The parameters of the model are also set to closely match \cite{lu2024boosting}.
Request data was obtained from CAT as detailed above, and time was discretized in minutes.
Service times for both pickups and drop-offs are set to 5 minutes, and road travel times are obtained from GraphHopper~\cite{GraphHopper}.
Time windows for request nodes are defined around the \emph{requested arrival time}: drop-offs are allowed from 30 minutes before up to the requested arrival time, and the pickup time window is the same interval shifted back in time by the direct driving time.
The time windows at the depot nodes are set from 5am-10pm to match the CAT service hours and the maximum shift length $L$ is set to 8 hours.
In experiments where the shifts are fixed, these fixed shifts are used as the depot time windows instead.
Shifts are restricted to start at the top of the hour, which is enforced through $\mathcal{T} = \{5\textrm{am}, 6\textrm{am}, \hdots, 2\textrm{pm}\}$.
Vehicle capacity is set to $Q=3$ passengers.
Finally, the large constant $M$ is set to $10,000$. This value strictly exceeds the total number of minutes in a day ($1,440$), guaranteeing a hierarchical objective structure that prioritizes fulfilling requests over minimizing time without introducing the unintended trade-offs of a weighted objective.

\section{Results}
\label{sec:results}

This section evaluates the new CP model on a real-world paratransit case study in Savannah.
Performance is compared both to the current practice at CAT and to the recently proposed algorithm by \cite{lu2024boosting}.
It also examines how the acceleration technique introduced in the methodology section benefits the search, and quantifies the value of allowing fully flexible driver shifts that are not restricted to start at the top of the hour.
The analysis will be presented through Tables~\ref{table:sota}-\ref{table:acceleration}, which all have a similar structure.
The \textit{Day} column denotes the calendar date of each instance; \textit{Vehicles} indicates the fleet size used by the \cite{lu2024boosting} baseline (matched to ensure a fair comparison); \textit{Requests} gives the total number of trip requests; \textit{Served} reports the number of requests fulfilled by each algorithm; and \textit{Time} records the CPU time in seconds used to obtain each solution.
Note that the CP model is configured to continue searching until the time limit, and thus always uses exactly 1,800 seconds. Because the OR-Tools routing solver is used here as a heuristic local-search procedure, it does not return a meaningful optimality gap for these runs; the reported served counts correspond to the best incumbent found within the time limit.

\begin{table}[p]
	\centering
	\caption{Performance of current practice (CAT), the algorithm by \cite{lu2024boosting}, and the CP model. The CP model either uses the shifts provided by \cite{lu2024boosting} or optimizes them directly. The ``Best'' row counts the number of instances (out of 15) on which a method attains the maximum number of served requests.}
	\begin{adjustbox}{width=\linewidth,center}
		\begin{tabular}{lccccccccc}
			\toprule
			& & & \multicolumn{1}{c}{CAT} & \multicolumn{2}{c}{\cite{lu2024boosting}} & \multicolumn{2}{c}{CP (shifts provided)} & \multicolumn{2}{c}{CP (shifts optimized)}\\
			\cmidrule(lr){4-4} \cmidrule(lr){5-6} \cmidrule(lr){7-8} \cmidrule(lr){9-10}
			Day & Vehicles & Requests & Served & Served & Time (s) & Served & Time (s) & Served & Time (s)\\
			\midrule
			20191202 & 30 & 475 & 365 & 437 & 1,534 & 440 & 1,800 & 418 & 1,800 \\
			20191203 & 29 & 462 & 342 & 439 &   898 & 447 & 1,800 & 432 & 1,800 \\
			20191204 & 35 & 561 & 450 & 519 & 2,526 & 525 & 1,800 & 510 & 1,800 \\
			20191205 & 32 & 515 & 409 & 480 & 1,737 & 478 & 1,800 & 466 & 1,800 \\
			20191206 & 32 & 509 & 391 & 460 & 1,241 & 456 & 1,800 & 446 & 1,800 \\
			20191209 & 31 & 493 & 427 & 461 & 1,594 & 465 & 1,800 & 444 & 1,800 \\
			20191210 & 32 & 508 & 419 & 478 & 1,571 & 492 & 1,800 & 484 & 1,800 \\
			20191211 & 34 & 549 & 443 & 508 & 1,871 & 508 & 1,800 & 511 & 1,800 \\
			20191212 & 32 & 506 & 421 & 478 & 1,824 & 483 & 1,800 & 462 & 1,800 \\
			20191213 & 32 & 511 & 398 & 479 & 1,231 & 479 & 1,800 & 463 & 1,800 \\
			20191216 & 30 & 478 & 405 & 440 & 1,195 & 435 & 1,800 & 444 & 1,800 \\
			20191217 & 29 & 456 & 405 & 435 &   986 & 428 & 1,800 & 422 & 1,800 \\
			20191218 & 30 & 474 & 384 & 442 & 1,491 & 438 & 1,800 & 434 & 1,800 \\
			20191219 & 32 & 518 & 443 & 483 & 1,415 & 489 & 1,800 & 482 & 1,800 \\
			20191220 & 31 & 497 & 363 & 456 & 1,397 & 469 & 1,800 & 469 & 1,800 \\
			\midrule
			\multicolumn{3}{r}{Average:}         &  404    & 466 & 1,501 & 463 & 1,800 & 459 & 1,800 \\
			\multicolumn{3}{r}{Percentage served:}&  80.73\%    & 93.12\% & -- & 92.41\% & -- & 91.68\% & -- \\
			\multicolumn{3}{r}{Best:}            &  0/15    & 5/15 & -- & 7/15 & -- & 3/15 & -- \\
			\bottomrule
		\end{tabular}
	\end{adjustbox}
	\label{table:sota}

\end{table}
\begin{table}[p]
	\centering
	\caption{Comparison between Original and Accelerated objective functions. The ``Best'' row uses the same convention as Table~\ref{table:sota}.}
	\label{table:acceleration}
	\begin{tabular}{lccccc}
		\toprule
		Day     & Vehicles & Requests & CAT & Original (Eq.~\ref{eq:model:obj}) & Accelerated (Eq.~\ref{eq:obj:better})\\
		\midrule
		20191202 & 30 & 475 & 365 & 388 & 418 \\
		20191203 & 29 & 462 & 342 & 389 & 432 \\
		20191204 & 35 & 561 & 450 & 483 & 510 \\
		20191205 & 32 & 515 & 409 & 437 & 466 \\
		20191206 & 32 & 509 & 391 & 421 & 446 \\
		20191209 & 31 & 493 & 427 & 415 & 444 \\
		20191210 & 32 & 508 & 419 & 470 & 484 \\
		20191211 & 34 & 549 & 443 & 470 & 511 \\
		20191212 & 32 & 506 & 421 & 440 & 462 \\
		20191213 & 32 & 511 & 398 & 432 & 463 \\
		20191216 & 30 & 478 & 405 & 395 & 444 \\
		20191217 & 29 & 456 & 405 & 396 & 422 \\
		20191218 & 30 & 474 & 384 & 401 & 434 \\
		20191219 & 32 & 518 & 443 & 451 & 482 \\
		20191220 & 31 & 497 & 363 & 436 & 469 \\
		\midrule
		\multicolumn{3}{r}{Average:}          & 404   & 428       & 459       \\
		\multicolumn{3}{r}{Percentage served:} & 80.73\% & 85.52\% & 91.68\% \\
		\multicolumn{3}{r}{Best:}              & 0/15  & 0/15      & 15/15     \\
		\bottomrule
	\end{tabular}
\end{table}

\subsection{Comparison to Current Practice and the State of the Art}
Table~\ref{table:sota} compares performance between four different methods: the current practice at Chatham Area Transit (CAT), the algorithm by \cite{lu2024boosting}, and two variants of the new CP model.
The first variant, \emph{CP (shifts provided)}, takes the shifts determined by \cite{lu2024boosting} and provides them to the CP model as a fixed input.
It only remains for the CP model to optimize the routes.
This setting is useful to see if the solution by \cite{lu2024boosting} can be improved further, and whether the CP model performs well in practical settings where the shifts have already been fixed.
The second variant, \emph{CP (shifts optimized)}, assumes no prior information and determines the shifts completely from scratch as part of the optimization.

% On the case study,
The CP model strictly outperforms current practice and is able to serve an average of 45 additional requests per day, increasing the overall service rate by over 10\%.
This result holds true whether the shifts are fixed and provided, or whether they are optimized on the fly.
The performance of the CP model is also remarkably close to that of the algorithm by \cite{lu2024boosting}.
Their advanced machine-learning boosted column-generation algorithm serves the most requests on average, but the relatively simple CP model is able to come within 1.5\% of this performance.

The last two columns of Table~\ref{table:sota} show that the CP model also delivers strong performance when used completely independently.
Fixing the shifts to those generated by \cite{lu2024boosting} does provide additional value and increases the average service rate by approximately 1\%.
The other way around, the comparison between \cite{lu2024boosting} and \emph{CP (shifts provided)} shows that the CP model may bring additional value to the algorithm by \cite{lu2024boosting}.
For day 20191210, for example, the CP model is able to take the shifts determined by \cite{lu2024boosting} and increase the number of requests served by 14 (3\% improvement).
This suggests that the CP model is able to find solutions that are hard to find with other methods.

\begin{figure}[!ht]
	\centering
	\includegraphics[width=\linewidth, trim={0 0 0 1cm}, clip]{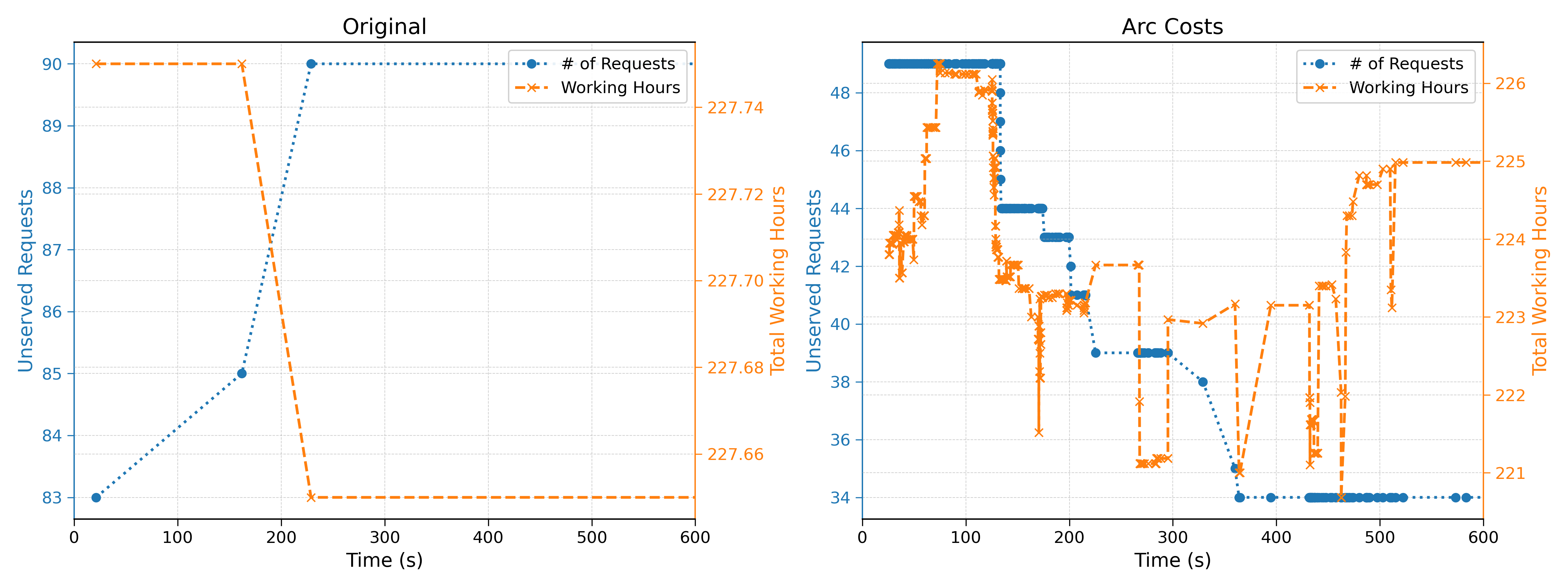}
	\caption{Requests served and total working hours during the solving process for instance 20191216. Original (left) versus Accelerated (right).}
	\label{fig:obj_curve_side_by_side}
\end{figure}

\subsection{Benefit of the Acceleration Technique}
This subsection provides an ablation study for the modeling choice introduced in Section~\ref{sec:acceleration}.
Table~\ref{table:acceleration} empirically demonstrates the importance of the acceleration technique to obtain the results above.
The \emph{Original} column provides the performance when Objective~\eqref{eq:model:obj} is used directly, while \emph{Accelerated} uses the proposed Equation~\eqref{eq:obj:better}, which includes total service time and travel time as a secondary objective.
There is a stark difference between the methods: \emph{Accelerated} outperforms \emph{Original} on every instance, and the average improvement is 31 additional requests served per day.
This corresponds to more than 6\% improvement in the overall service rate.
It is clear that the augmented objective is important for the CP model to compete with the approach by \cite{lu2024boosting}.

Figure~\ref{fig:obj_curve_side_by_side} provides insight into how the acceleration technique impacts the behavior of the solver for a representative instance.
The plots show the number of unserved requests and the total working hours against elapsed CPU time as the solver progresses. The reported curves correspond to the current accepted search solution rather than the best incumbent retained, so the unserved-request count may temporarily increase when tabu search accepts a worsening move; the best incumbent is monotonically non-increasing.
Total working hours are computed as the sum of travel times between depots for all active vehicles and may increase when new requests are inserted or when additional vehicles are deployed.
Under Objective~\eqref{eq:model:obj} (left panel), the solver constructs an initial solution with 83 unserved requests, but is unable to find better solutions after that.
In contrast, the augmented objective (right panel) not only produces a significantly better initial solution, but also results in a steady decline in unserved requests as the solver progresses.
Note that there are several occasions where the number of requests plateaus.
However, progress is not stalled, as the solver focuses on decreasing the number of working hours, which eventually creates the necessary space in the schedule to insert additional requests.

\begin{figure}[!ht]
    \centering
    \includegraphics[width=0.8\linewidth]{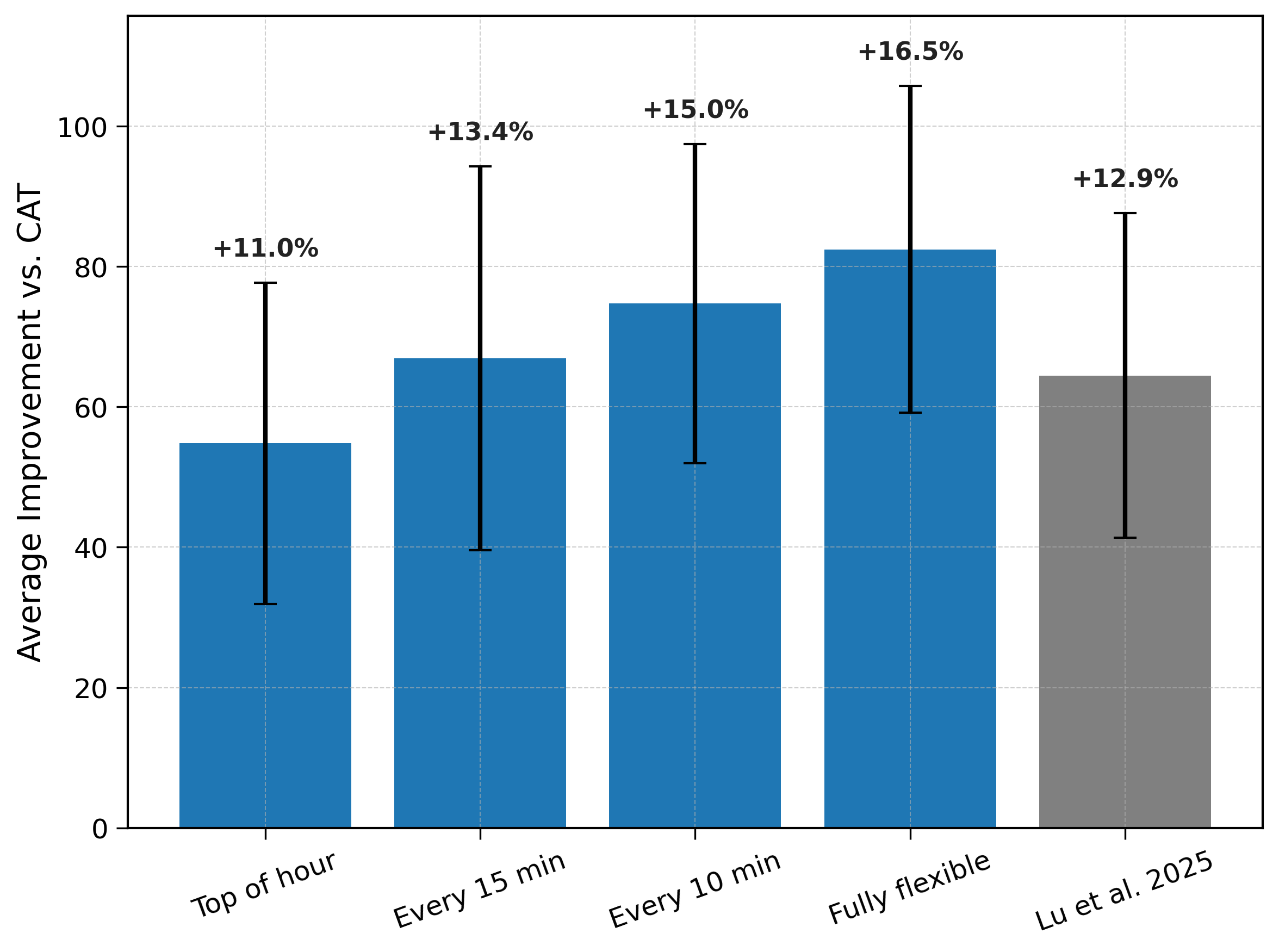}
    \caption{Average additional requests served compared to the CAT baseline, with percentage-point improvements over CAT's served percentage annotated for each CP flexibility level and \cite{lu2024boosting}.}
    \label{fig:flex}
\end{figure}

\subsection{Supporting More Flexible Shifts}
The new CP model also offers the opportunity to quantify the value of allowing more flexible driver shifts that are not constrained to start at the top of the hour.
This flexibility enables a closer alignment between driver availability and trip demand.
Constraints~\eqref{eq:model:start_times} restrict the set of possible driver shift start times to $\mathcal{T}$, which is currently defined as $\mathcal{T} = \{5\textrm{am}, 6\textrm{am}, \hdots, 2\textrm{pm}\}$ to enforce shifts starting exactly on the hour.
Increasing flexibility simply requires redefining the parameter $\mathcal{T}$ to include smaller intervals, for example $\mathcal{T} = \{5\textrm{am}, 5.15\textrm{am}, 5.30\textrm{am}, \hdots, 1.45\textrm{pm}, 2\textrm{pm}\}$ for 15-minute intervals, or even $\mathcal{T} = [5\textrm{am}, 2\textrm{pm}]$ for fully flexible shifts.

The ability to easily adjust the flexibility of driver shift start times is a major advantage of CP compared to the column-generation approach by \cite{lu2024boosting}, which requires solving a separate subproblem for each possible start time.
This means that computation time grows with the cardinality of $\mathcal{T}$, making highly flexible settings increasingly difficult to manage.
A fully flexible setting could potentially be addressed with time-dependent cost functions (e.g., see \cite{IoachimEtAl1998-DynamicProgrammingAlgorithm}), but this would further complicate the algorithm and likely increase computation time significantly, whereas increasing flexibility in the CP framework is almost trivial.

Figure~\ref{fig:flex} shows that increasing shift-start flexibility leads to clear and substantial improvements in operational performance. Allowing shifts to begin every 15 minutes or every 10 minutes raises the average number of requests served by 67 and 75 per day, respectively, compared to the CAT baseline. These intermediate levels of flexibility already eliminate a large share of the unmet demand that remains under the conventional top-of-the-hour shift structure used in both the current system and in \cite{lu2024boosting}. Full flexibility delivers the greatest benefit: the fully flexible configuration improves average performance by 83 additional requests per day, corresponding to a 16.5\% increase in service rate, achieving near-complete demand fulfillment. Moreover, under full flexibility, the CP model achieves a substantially higher average number of served requests than the approach proposed in \cite{lu2024boosting}.

While the gains from fully flexible driver shifts can be substantial, they should be balanced against the need for more sophisticated driver timekeeping systems and greater driver adaptability.
The new CP model can be used to explore these trade-offs, and can help stakeholders to assess the value of shift policies ranging from completely fixed to fully flexible, and everything in between.

\subsection{Sensitivity Analysis} \label{sec:sensitivity}
This subsection provides empirical support for the implementation choices in Section~\ref{sec:methodology} by varying the local-search metaheuristic and the first-solution strategy used by the routing solver, with all other settings fixed to those of the \emph{CP (shifts optimized)} configuration.

Table~\ref{tab:metaheuristic-sensitivity} compares local-search metaheuristics. The standard Tabu Search achieves the highest average number of served requests (459.87), while Guided Local Search yields the lowest average total working hours (237.41). The default \texttt{GENERIC\_TABU\_SEARCH} provides a balanced configuration with 459.13 served requests and 237.50 working hours, and is used throughout the rest of the paper.

Table~\ref{tab:first-solution-sensitivity} compares first-solution strategies. The Automatic and Parallel Cheapest Insertion strategies are consistently strongest, reaching an average of 451.87 initial and 459.13 final served requests across all 15 instances. By contrast, Path Cheapest Arc fails to produce any feasible solution, and Savings, All Unperformed, Global Cheapest Arc, and First Unbound Min Value drop the final service level to 30 or fewer requests on average. These results justify retaining the Automatic default.

\begin{table}[!ht]
\centering
\caption{Metaheuristic comparison for CP (shifts optimized). Reported values are averages across all available dates. The default method is bolded.}
\label{tab:metaheuristic-sensitivity}
\begin{tabular}{lrr}
\toprule
Metaheuristic & Avg.\ served & Avg.\ total working hours \\
\midrule
Automatic & 459.13 & 237.60 \\
Greedy Descent & 459.13 & 237.60 \\
Guided Local Search & 459.27 & 237.41 \\
Simulated Annealing & 459.20 & 237.60 \\
Tabu Search & 459.87 & 237.98 \\
\textbf{Generic Tabu Search} & 459.13 & 237.50 \\
\bottomrule
\end{tabular}
\end{table}

\begin{table}[!ht]
\centering
\caption{First-solution strategy comparison for CP (shifts optimized). Reported values are averages across all available dates. The default method is bolded.}
\label{tab:first-solution-sensitivity}
\begin{adjustbox}{width=\linewidth,center}
\begin{tabular}{lrrrr}
\toprule
First-solution strategy & Avg.\ initial served & Avg.\ time to first feasible (s) & Avg.\ final served & Feasible runs \\
\midrule
\textbf{Automatic} & 451.87 & 32.23 & 459.13 & 15 \\
Path Cheapest Arc & -- & -- & -- & 0 \\
Savings & 1.71 & 22.73 & 29.93 & 14 \\
All Unperformed & 0.00 & 0.27 & 29.00 & 15 \\
Parallel Cheapest Insertion & 451.87 & 31.72 & 459.13 & 15 \\
Local Cheapest Insertion & 420.27 & 35.63 & 432.53 & 15 \\
Global Cheapest Arc & 0.00 & 2.09 & 29.00 & 15 \\
First Unbound Min Value & 0.00 & 1.34 & 29.00 & 15 \\
\bottomrule
\end{tabular}
\end{adjustbox}
\end{table}

\section{Conclusion}
\label{sec:conclusion}

This paper considered the problem of jointly optimizing route planning and shift scheduling for paratransit services.
These services are vital for individuals who cannot use fixed-route public transit, including those with disabilities.
Optimizing these services is therefore essential for transit agencies to deliver high-quality service efficiently.

A CP model was introduced to solve this problem, along with practical guidance to implement the model in Google OR-Tools for real-world use.
This complements prior work by \cite{lu2024boosting}, which presented a highly effective but technically complex combination of AI and OR that also relies on historical data.
The new model is significantly easier to implement and provides an inherently practical solution for transportation planners.

Through a case study in Savannah, it was demonstrated that the new approach is competitive with the approach by \cite{lu2024boosting} and is able to come within 1.5\% of the performance of the more complex method.
Both methods increase the overall service rate of current practice by over 10\%.
The CP model was shown to be effective both whether driver shifts were provided or determined as part of the optimization.
It was also found that the CP model may bring value even if the method by \cite{lu2024boosting} is already implemented: running the model as a post-processing step may still increase the number of requests served by up to 3\%.

Analysis showed that the introduced acceleration technique is critically important to compete with the approach by \cite{lu2024boosting}.
Combining the number of requests and the total working hours in a single objective helps guide the solver towards more time-efficient routes, which in turn allows more requests to be inserted.
This finding provides empirical support for the modeling choice in Section~\ref{sec:acceleration}, rather than treating it as a solver-independent best practice.

Finally, the paper quantified the value of allowing fully flexible drivers shifts that do not necessarily start at the top of the hour, enabling closer alignment between driver availability and trip demand.
It is a feature of the CP model that this can be achieved by simply removing some of the constraints.
This added flexibility yields an additional 5\% improvement in the number of requests served.
It also allows the CP model to outperform the algorithm by \cite{lu2024boosting}, for which it is not trivial to offer this kind of flexibility.

Future research in this field could focus on further improving the models to incorporate additional practical constraints and further improve solution quality.
Another interesting direction would be the implementation of a more dynamic system that reoptimizes the plan as new requests come in and as things change during the day.
This dynamic piece would help account for randomness due to traffic, cancellations, and changes to rider itineraries.
The CP model presented in this paper may provide a basis for such a system.

\section*{Acknowledgements}
This research is partly supported by NSF awards 2112533 and 1854684.

\bibliography{references}

\end{document}